\documentclass{amsart}
\usepackage{amsmath}
\usepackage{amssymb}
\usepackage{amscd}
\newtheorem{theorem}{Theorem}[section]
\newtheorem{lemma}[theorem]{Lemma}
\newtheorem{proposition}[theorem]{Proposition}
\newtheorem{corollary}[theorem]{Corollary}

\newtheorem*{addendum}{Theorem 6.1 Addendum}

\theoremstyle{definition}
\newtheorem{definition}[theorem]{Definition}

\newtheorem{example}[theorem]{Example}

\newtheorem{definitions and remarks}[theorem]{Definitions and Remarks}

\theoremstyle{remark}
\newtheorem{remark}[theorem]{Remark}
\newtheorem{remarks}[theorem]{Remarks}

\numberwithin{equation}{section}

\newcommand{\inv}{\mathrm{inv}}

\newcommand{\oa}{\overline{a}}

\newcommand{\Spec}{\mathrm{Spec}\,}

\newcommand{\Sing}{\mathrm{Sing}\,}
\newcommand{\supp}{\mathrm{supp}\,}
\newcommand{\cosupp}{\mathrm{cosupp}\,}

\newcommand{\ord}{\mathrm{ord}}

\newcommand{\car}{\mathrm{char}\,}

\newcommand{\length}{\mathrm{length}\,}

\newcommand{\Der}{\mathrm{Der}}

\newcommand{\Bl}{\mathrm{Bl}}
\newcommand{\loc}{\mathrm{loc}}

\newcommand{\al}{{\alpha}}

\newcommand{\g}{{\gamma}}

\newcommand{\ka}{{\kappa}}

\newcommand{\La}{{\Lambda}}

\newcommand{\p}{{\partial}}
\newcommand{\s}{{\sigma}}

\newcommand{\vp}{{\varphi}}

\newcommand{\IN}{{\mathbb N}}
\newcommand{\IQ}{{\mathbb Q}}
\newcommand{\IA}{{\mathbb A}}

\newcommand{\IC}{{\mathbb C}}

\newcommand{\IF}{{\mathbb F}}

\newcommand{\cC}{{\mathcal C}}
\newcommand{\cD}{{\mathcal D}}

\newcommand{\cF}{{\mathcal F}}
\newcommand{\cG}{{\mathcal G}}

\newcommand{\cI}{{\mathcal I}}
\newcommand{\cJ}{{\mathcal J}}

\newcommand{\cO}{{\mathcal O}}

\newcommand{\uk}{\underline{k}}
\newcommand{\ul}{\underline{l}}
\newcommand{\um}{\underline{m}}
\newcommand{\fm}{{\mathfrak m}}
\newcommand{\fp}{{\mathfrak p}}
\newcommand{\fC}{{\mathfrak C}}
\newcommand{\fN}{{\mathfrak N}}
\newcommand{\fV}{{\mathfrak V}}

\newcommand{\wcO}{{\widehat \cO}}

\newcommand{\wX}{{\widehat X}}
\newcommand{\wY}{{\widehat Y}}

\newcommand{\ucG}{\underline{\cG}}

\newcommand{\ucI}{\underline{\cI}}
\newcommand{\ucJ}{\underline{\cJ}}
\newcommand{\ucC}{\underline{\cC}}
\newcommand{\ucD}{\underline{\cD}}

\newcommand{\uX}{\underline{X}}

\newcommand{\lbr}{{[\![}}
\newcommand{\rbr}{{]\!]}}

\begin{document}

\title[$\IQ$-universal desingularization]
{$\IQ$-universal desingularization}

\author{Edward Bierstone}
\address{Department of Mathematics, University of Toronto, 40 St. George Street,
Toronto, Ontario, Canada M5S 2E4}
\email{bierston@math.toronto.edu}
\thanks{The first two authors' research was supported in part by NSERC
grants OGP0009070 and OGP0008949.}

\author{Pierre D. Milman}
\address{Department of Mathematics, University of Toronto, 40 St. George Street,
Toronto, Ontario, Canada M5S 2E4}
\email{milman@math.toronto.edu}

\author{Michael Temkin}
\address{Department of Mathematics, University of Pennsylvania, Philadelphia,
Pennsylvania 19104, USA}
\email{temkin@math.upenn.edu}



\begin{abstract}
We prove that the algorithm for desingularization of algebraic varieties
in characteristic zero of the first two authors is functorial with respect
to regular morphisms. For this purpose, we show that, in characteristic zero, 
a regular morphism with connected affine source
can be factored into a smooth morphism, a ground-field extension and a
generic-fibre embedding. Every variety of characteristic zero admits a 
regular morphism to a $\IQ$-variety. The desingularization algorithm is therefore
\emph{$\IQ$-universal} or \emph{absolute} in the sense that it is induced from
its restriction to varieties over $\IQ$. As a consequence, 
for example, the algorithm
extends functorially to localizations and Henselizations of varieties. 
\end{abstract}

\maketitle
\setcounter{tocdepth}{1}
\tableofcontents

\section{Introduction}
Our main result is the following.

\begin{theorem}
Every algebraic variety in characteristic zero admits
(strong) resolution of singularities that is functorial with respect
to regular morphisms.
\end{theorem}

More precisely, we show that the
desingularization algorithm of \cite{BMinv, BMfunct} is functorial
with respect to regular morphisms. (See Theorem and Addendum 6.1 below
for a precise statement of Theorem 1.1. ``Strong'' means in particular
that the desingularization is by blowings-up along smooth subvarieties.)
The assertion of Theorem 1.1 is called
$\IQ$\emph{-universal} resolution of singularities
by Hironaka \cite{Hid} because
any algebraic variety $X$ in characteristic zero admits a regular
morphism to a variety $Y$ defined over the rational numbers $\IQ$
(see Theorem 3.1 below), so that resolution of singularities of $X$ is
induced by that of $Y$. In \cite{Hid}, Hironaka writes that $\IQ$-universal
desingularization will be proved in a subsequent paper, but a proof
has not appeared before as far as we know (see Remark 5.5).

An \emph{(algebraic) variety} means a scheme $X$ which admits a morphism
of finite type $X \to \Spec \uk$, where $\uk$ is a field. (It will be
convenient to extend this definition to schemes that are disjoint unions
of such; see Remarks 4.7(2).) If a morphism $X \to \Spec \uk$ is fixed, 
we will say that $X$ is a variety with
\emph{ground field} $\uk$, or a $\uk$\emph{-variety}.

A morphism of schemes $f: X \to Y$ is \emph{regular} if $f$ is flat and
all fibres of $f$ are \emph{geometrically regular}; equivalently, if
$f$ is flat and, for every morphism $T \to Y$ of finite type, all fibres
of $X \times_Y T \to T$ are regular \cite[\S\S28,\,33]{Mat}. If $f$ is of
finite type, then $f$ is regular if and only if it is smooth
\cite[Thm.\,10.2]{Hart}. Thus regularity
is a generalization of smoothness to morphisms that are not necessarily of
finite type.

\begin{theorem}
A regular morphism $f: X \to Y$, where $X$ is a connected
affine variety and $Y$ is a variety over a field $\uk$ of characteristic
zero, can be factored as
\begin{equation}
X \,\cong\, Z_{\eta} \times_{\Spec \um} \Spec \ul \,\xrightarrow{f_{\ul}}\, Z_\eta
\,\xrightarrow{f_{\um}}\, Z \,\xrightarrow{f_{\uk}}\, Y,
\end{equation}
where $f_{\uk}$ is a smooth morphism of $\uk$-varieties,
$f_{\ul}$ is a ground-field
extension and $f_{\um}: Z_{\eta} \to Z$ is a \emph{generic-fibre embedding}.
\end{theorem}

A \emph{generic-fibre embedding} $f_{\um}: Z_{\eta} \to Z$
means there is a dominant $\uk$-morphism $Z \to T$ to an integral
$\uk$-variety $T$, and $f_{\um}: Z_{\eta} \to Z$ is the canonical morphism from
the \emph{generic fibre} $Z_{\eta} = Z \times_T \eta$ (where $\eta = \Spec \um$
is the generic point of $T$. See Sections 2 and 4.)

For example, $\Spec \IQ(x)[y] \to \Spec \IQ[x,y]$ is a generic-fibre embedding;
it is a regular morphism that is not a composite of smooth morphisms and
ground-field extensions.

Functoriality with respect to smooth morphisms and ground field extensions in
the strong desingularization algorithm for varieties
\cite{BMinv} is proved in \cite{BMfunct}. (\cite{W} and \cite{Ko} provide
versions of weak desingularization
of varieties that are also functorial with respect to smooth morphisms and
ground field extensions.) In Section 6, we deduce Theorem 1.1 from the
previous results, using Theorem 1.2
and functoriality with respect to generic-fibre embeddings (see \S4.3
and Proposition 6.3).

Note, however, that all previous results on functoriality seem to make a
tacit assumption that the smooth morphisms have constant relative dimension.
We impose no such restriction in Theorem 1.1, so we also have to show that
the desingularization algorithms of \cite{BMinv, BMfunct} are functorial with
respect to arbitrary smooth morphisms. We are grateful to Ofer Gabber for
raising this issue; see \S6.3.

Functorial desingularization involves important local--global issues.
For example, even if a variety has several connected
components (so that resolutions of singularities of different components are independent), functoriality depends on the order in which the
components are blown up. Such issues intervene throughout the article
(see \S4.3 and Remarks 5.1, 6.2).

The algorithm for strong resolution of singularities of \cite{BMinv, BMfunct}
is based on a desingularization algorithm for a \emph{marked ideal} (as presented
in \cite{BMfunct}). The proofs of the theorems involve a notion of
\emph{equivalence} of marked
ideals. (The meanings of these notions are recalled in \S5 below;
for details we refer to \cite{BMfunct}.)

\begin{theorem}
The algorithm for resolution of singularities of marked ideals in characteristic
zero (of \cite{BMfunct}) is functorial with respect to equivalence classes
(of marked ideals of a given dimension; see \S5.1) and with respect to
regular morphisms.
\end{theorem}

In Section 5, we obtain Theorem 1.3
from previous functoriality results (with respect to
smooth morphisms and ground field extensions) again using Theorem 1.2,
\S4.3 and functoriality with respect to generic-fibre embeddings (Propostion
5.3). If $Z_\eta \to Z$ is a generic-fibre embedding in characteristic zero,
where $Z$ is smooth, then equivalent marked ideals on $Z$
pull back to equivalent marked ideals on $Z_\eta$ (see Lemma 5.2).

An algorithm for
principalization of an ideal that is functorial with respect to regular
morphisms also follows from Theorem 1.3.

For Proposition 5.3, we follow the proof in \cite{BMfunct} step-by-step. The
only point that is not immediate involves passage
from a marked ideal to a (local) \emph{coefficient
ideal} (Step I in \cite{BMfunct}). Suppose that $\psi: Z_\eta \to Z$ is a
generic-fibre embedding as above, where $Z_\eta,\, Z$ are smooth. Let $\ucI_\eta$
denote the pullback $\psi^*(\ucI)$ to $Z_\eta$ of a marked ideal $\ucI$ on $Z$.
It follows from Lemma 5.1 that the coefficient ideal
of $\ucI$ pulls back to a marked ideal which is
equivalent to $\ucI_\eta$ (see Lemma 5.6).

\begin{remark}
To explain the significance of the latter, let us recall that the local coefficient
ideal for $\ucI$ is defined using ideals of derivatives of $\ucI$. If $\cI \subset
\cO_X$ is a coherent ideal on a regular variety $X$, then the \emph{derivative
ideal} $\cD(\cI)$ is generated by all
first derivatives of local sections of $\cI$. If $X$ is an $\um$-variety,
this means that $\cD(\cI)$ is the image of the natural morphism
$\Der_X \times \cI \to \cO_X$, where $\Der_X$ denotes the sheaf of
$\um$-derivations $\Der_{\um}(\cO_X,\cO_X)$;
i.e., $\um$-linear homomorphisms $\cO_X \to \cO_X$ that satisfy Leibniz's
rule (hence vanish on $\um$).

A coefficient ideal of a marked ideal on $Z$ involves $\uk$-derivations, while
on $Z_\eta$ a coefficient ideal involves $\um$-derivations. Derivative ideals
defined using $\uk$-derivations (or $\IQ$-derivations) may be much larger than
those defined over $\um$ because they involve derivatives along
``constants'' (elements of $\um$ that are transcendental over $\uk$ or $\IQ$).
\emph{$\IQ$-universal resolution of singularities} means that the derivative
ideals defined using $\IQ$-derivations nevertheless do not result in smaller
centres of blowing up, so we can run the
desingularization algorithm for a variety defined over a field or characteristic
zero, in general, using derivatives defined over $\IQ$.
\end{remark}

In Section 7, we present an alternative (though less explicit) approach to
universal desingularization algorithms based on approximation methods
of \cite[\S8]{EGAIV} that we use in our
proof of the factorization theorem 1.2. We can start with any desingularization
algorithm for varieties over $\IQ$ that is functorial with respect to smooth
morphisms, and extend it to a class of schemes over $\IQ$ that includes
all varieties of characteristic zero as well as their localizations and Henselizations
along closed subvarieties. The resulting
desingularization algorithm is again functorial with respect to regular morphisms.

\section{The generic fibre}
Let $\pi: Z \to T$ denote a dominant morphism of $\uk$-varieties, where
$\uk$ is a field and $T$ is integral. Let $\eta$ denote the generic point
of $T$; i.e., $\eta = \Spec \um$, where $\um$ is the function field $K(T)$ of $T$.
There is a fibred-product diagram
$$
\begin{CD}
Z \times_T \eta @>\psi >> Z \\
@VVV                   @VV\pi V \\
\eta    @>>>  T\\
\end{CD}
$$
in which all morphisms are dominant. The \emph{generic fibre} $Z_\eta$ of $\pi$
denotes the $\um$-variety $Z \times_T \eta$.

Suppose that $Z$ and $T$ are affine. Then $\um$ is the field of fractions
$\uk(T)$ of the coordinate ring $\uk[T]$; by definition, $\uk(T)$
is the localization
$\uk[T]_S$, where $S=S(T)$ is the multiplicative subset $\uk[T]\backslash \{0\}$
of $\uk[T]$. The morphism $\pi$ induces an injection $\uk[T] \hookrightarrow
\uk[Z]$, so that the coordinate ring of $Z_{\eta}$,
$$
\um[Z_{\eta}] = \uk(T) \otimes_{\uk[T]} \uk[Z] \cong \uk[Z]_S.
$$

\begin{remark} We recall that if $A$ is a ring and $S$ is a multiplicative
subset of $A$, then an ideal of $A_S$ is prime if and only if it of the
form $\fp\cdot A_S$, where $\fp$ is a prime ideal of $A$ disjoint from $S$.
\end{remark}

If $a \in Z$, we write $\fm_{Z,a}$ for the maximal ideal of $\cO_{Z,a}$ and
$\ka(a)$ for the residue field $\cO_{Z,a}/\fm_{Z,a}$. Let
$b$ be a point of $Z_\eta$ and let $a = \psi(b) \in Z$. ($a$ is not necessarily
closed, even if $b$ is closed.) Let $\psi^*_b: \cO_{Z,a} \to \cO_{Z_{\eta},b}$
denote the homomorphism of local rings induced by $\psi$.

If $T$ is affine as above, then $\cO_{Z_{\eta},b}$ can be 
identified with the localization $(\cO_{Z,a})_S$, where $S=S(T)$;
thus $\fm_{Z,a} \cap S = \emptyset$, by Remark 2.1, so that
$\cO_{Z_{\eta},b} \cong (\cO_{Z,a})_S =\cO_{Z,a}$. Therefore,
(1) $\ka(a) = \ka(b)$;  (2) $\cO_{Z_{\eta},b}$ is
regular if and only if $\cO_{Z,a}$ is regular; (3) if $f \in \cO_{Z,a}$, then
the order $\ord\, f = \ord\, \psi^*_b(f)$.

Suppose that $W$ is an integral subvariety of $Z$. It follows from
Remark 2.1 that $\pi|_W: W \to T$ is dominant if and only if there is a
point $b \in Z_\eta$ such that $W$ is the closure $\oa$ of $a = \psi(b)$.

Suppose that $b$ is a closed point of $Z_\eta$. Then $\ka(b)/\um$ is a
finite field extension, by the Nullstellensatz \cite[Th. 4.19]{Eisen}.
Let $a = \psi(b)$ and let $W=\oa$. Then $\ka(b) = \ka(a) = K(W)$, so that
$K(W)/K(T)$ is a finite extension, and $\pi|_W$ is a 
generically finite morphism.

In general, if $\pi|_W$ is not dominant, then $\psi^{-1}(W) = \emptyset$,
and if $\pi|_W$ is dominant, then $\psi^{-1}(W) = W_\eta$. In the latter
case, if $\uk$ is perfect, then $W_\eta$ is smooth if and only if there
is an open subset of $T$ over which $W$ is smooth (and the restriction of $\pi$
is a smooth morphism).

\section{Embedding as the generic fibre of a $\IQ$-variety}
Every variety $X$ in characteristic zero can be obtained by a
base change from a variety which admits a generic fibre embedding
$Z_{\eta} \to Z$ into a variety $Z$ over $\IQ$. This is a well-known result
which is a special case of Theorem 1.2. We outline a proof for completeness and
also to illustrate a technique
developed in great generality by Grothendieck \cite[Th.\,8.8.2, Prop.\,8.13.1]
{EGAIV} that we will use to prove Theorem 1.2.

\begin{theorem}
Let $X$ denote a variety over a field $\uk$ of characteristic zero.
Then there exists a $\IQ$-variety $Z$ and a dominant morphism $\pi: Z \to T$,
where $T$ is an integral $\IQ$-variety, such that
$X$ is obtained from the generic fibre $Z_\eta$ of
$\pi$ by base extension; i.e., $X = Z_{\eta} \times_{\Spec \um} \Spec \uk$,
where $\um = K(T)$.
If $X$ is smooth, then we can take $Z$ smooth.
\end{theorem}

\begin{proof}
Our $\uk$-variety $X$ can be constructed by glueing together finitely many
affine varieties $X_i$ along open subsets. (See \cite[Ch.\,II, Exercise 2.12]
{Hart}.) Let $\um$ denote the subfield of $\uk$ obtained by extending $\IQ$
by the coefficients of the polynomials comprising (finite) generating sets
for the ideals $I_i$, where $X_i = \Spec \uk[y]/I_i$, for all i, together
with the coefficients of the polynomials needed to present the glueing data.
In other words: Let $\{c_j\} \subset \uk$ denote the (finite) set of all
coefficients above, and consider the ring homomorphism $\g: \IQ[x] \to \uk$
given by $\g(x_j) = c_j$, where $\IQ[x]$ denotes the ring of polynomials over
$\IQ$ in indeterminates $x_j$. The $\ker \g$ is a prime ideal $\fp$, and
$\um$ denotes the field of fractions of $\IQ[x]/\fp$.

The field $\um$ is an extension of $\IQ$ of finite type. Our variety $X$
can be considered also as a variety $Z_{\um}$ defined over $\um$. (As a
$\uk$-variety, $X$ is obtained from $Z_{\um}$ by base extension $\Spec \uk
\to \Spec \um$.)

Let $T := \Spec \IQ[x]/\fp$. For each affine chart $X_i = \Spec \uk[y]/I_i$
above, let $J_i \subset \IQ[x,y]$ denote the ideal with generators obtained
from those above by replacing each coefficient $c_j$ by $x_j$. Then there is
a $\IQ$-variety $Z$ constructed by glueing together the affine varieties
$Z_i = \Spec \IQ[x,y]/(J_i, \fp)$ (where $(J_i, \fp)$ denotes the ideal generated
by $J_i$ and $\fp$) using glueing morphisms obtained in the same way from those
for $X = \cup X_i$.

Clearly, there is a dominant morphism $\pi: Z \to T$, and $Z_{\um}$ can be
identified with the generic fibre $Z_\eta$ of $\pi$ (where $\eta = \Spec \um$
is the generic point of $T$.)

If $X$ is smooth, then $Z_{\um}$ is smooth. The variety $Z$ is
\emph{a priori} singular, but we can restrict to an open subset of $T$ over
which it is smooth (as in \S2).
\end{proof}

\section{Factorization of a regular morphism}

\subsection{Ground fields}
A variety $X$ may admit many different structures of a $\uk$-variety (i.e.,
many different morphisms of finite type $f:X\to\Spec \uk$, even for a fixed
field $\uk$). This is usually the case when $X$ is
not reduced. (A simpler possibility is that $\uk$ is finite over a subfield
isomorphic to $\uk$ itself).
Nevertheless, a connected reduced variety possesses a  unique maximal
ground field, so there is a natural choice of ground field.

\begin{lemma}\label{conredlem}
Let $X$ be a connected reduced variety. Then the ring $\cO_X(X)$ contains a
maximal subfield $\uk$ containing
any other subfield. In particular, any morphism $X\to\Spec \ul$, where $\ul$ is
a field, factors through the morphism
$X\to\Spec \uk$ corresponding to the embedding $\uk\hookrightarrow\cO_X(X)$.
\end{lemma}

\begin{proof}
Let $\IF$ be the prime subfield contained in $\cO:=\cO_X(X)$ and let $\uk$ be the set of elements $f\in\cO$ such
that $\cO$ contains the subfield $\IF(f)$. It suffices to prove that $\uk$ is a subfield of $\cO$, since it is
then clear that $\uk$ is as required. Fix a structure $X\to\Spec \ul$ of an $\ul$-variety on $X$. Then any
element $f\in\cO$ induces a morphism $F:X\to\IA^1_{\ul} = \Spec \ul[T]$ whose image $Z$ is constructible, by
Chevalley's theorem \cite[Ex.\,3.19]{Hart}. Since $Z$ is connected, either it is a point or it omits at most
finitely many points. In the latter case, $f\notin \uk$ because $\IF[T]$ has infinitely many primes. On the
other hand, in the first case, the map $F$ is constant on $X$ and equal
to %
an element of the algebraic closure of $\ul$, hence $f$ annihilates an irreducible polynomial over $\ul$ and so
$\ul[f]$ is a subfield of $\cO$. This proves that $F$ is constant if and only if $f\in \uk$. It follows that
$\uk$ is a ring, and we have also seen that $\ul[f]$ is a field for any element $f\in \uk$. Thus $\uk$ is a
field as required. (We have actually proved that it coincides with the integral closure of $\ul$ in $\cO$.)
\end{proof}

The lemma following shows that to a general connected affine variety $X$,
we can assign a uniquely defined
maximal ground field $\uk$ just by taking the maximal ground field of its
reduction. However, in sharp contrast to
the reduced case above, the structure morphism $X\to\Spec \uk$ is absolutely
not canonical. In particular,
non-isomorphic $\uk$-varieties can be isomorphic as abstract schemes.

\begin{lemma}\label{conlem}
Let $X$ be an affine variety with reduction $X_0$. Then any morphism
$f_0:X_0\to\Spec \uk$, where $\uk$ is a field, extends
to a morphism $f:X\to\Spec \uk$ in the sense that $f_0$ is the composition
of $f$ with the reduction morphism
$X_0\to X$. Moreover, suppose that $\ul$ is a subfield of $\uk$ such that
$\uk/\ul$ is separable, and fix an extension $f':X\to\Spec \ul$
of the morphism $f'_0:X_0\to\Spec \uk \to \Spec \ul$. Then we can choose
$f$ compatible with $f'$.
\end{lemma}

\begin{proof}
Let $\IF$ be the prime subfield of $\uk$. Then $\uk/\IF$ is separable and $X$ admits a unique morphism to $\Spec
\IF$. Let $\ul$ be a subfield of $\uk$ such that $\uk/\ul$ is separable (for example, $\ul = \IF$ to get the
first assertion of the lemma). Since $\uk/\ul$ is separable, the morphism
$\Spec \uk\to\Spec \ul$ 
is formally smooth, by \cite[Prop.\,28.I]{Mat}. In particular, the $\ul$-morphism $X_0\to\Spec \uk$ extends
to %
an $\ul$-morphism $X\to\Spec \uk$.
\end{proof}

Combining Lemmas \ref{conredlem} and \ref{conlem},
we obtain the following corollary.

\begin{corollary}\label{concor}
Let $X$ be a connected affine variety. Then any morphism $X\to\Spec \ul$, where
$\ul$ is a perfect field, factors through
a morphism $X\to\Spec \uk$ of finite type, where $\uk$ is a field.
\end{corollary}

The following example shows why we have to assume that $X$ is connected
in Corollary \ref{concor}.

\begin{example}
Let $\uk$ be a field which admits an endomorphism $\vp: \uk \to \uk$ with
$[\uk : \vp(\uk)] = \infty$; for example, $\uk = \IC$. Set $S := \Spec \uk$ and
let $X = S_1 \coprod S_2$ be the disjoint union of two copies of $S$. Then the
morphism $X \to S$ which restricts to the identity on $S_1$ and to $\Spec \vp$
on $S_2$ cannot be factored as in Corollary \ref{concor}.
\end{example}

\subsection{Regular morphisms}
The assertion of our factorization theorem 1.2 is included in the
following result.

\begin{theorem}\label{regprop}
Let $\uk$ denote a perfect field, and let $f: X \to Y$ denote a morphism, where
$X$ is a connected affine variety and $Y$ is a $\uk$-variety.
\begin{enumerate}
\item
The $f$ can be factored as
\begin{equation}
X \,\cong\, Z_{\eta} \times_{\Spec \um} \Spec \ul \,\xrightarrow{f_{\ul}}\, Z_\eta
\,\xrightarrow{f_{\um}}\, Z \,\xrightarrow{f_{\uk}}\, Y,
\end{equation}
where $f_{\uk}$ is a morphism of $\uk$-varieties,
$f_{\ul}$ is a ground-field
extension and $f_{\um}$ is a generic-fibre embedding.
\item
Assume that $\car \uk = 0$. Then $f$ is regular if and only if the morphism
$f_{\uk}$ in (4.1) is
smooth on a neighbourhood $U$ of $Z_\eta$. (So, if $f$ is regular, then
we get (4.1) with $f_k$ smooth by restricting to $U$.) 
\end{enumerate}
\end{theorem}

\begin{proof}
By Corollary \ref{concor} the morphism $X \xrightarrow{f}  Y \to \Spec \uk$
extends to a morphism $X \to \Spec \ul$ of finite type. We
construct $Z_\eta$ by approximating $X$ with a variety defined over a
finitely generated $\uk$-field $\um$.
The $\uk$-scheme $\Spec\ul$ is the projective limit of the $\uk$-schemes
$\Spec\um_i$ where $\um_i$ runs over all
subfields of $\ul$ that contain $\uk$ and are finitely generated over $\uk$.
By \cite[Thm.\,8.8.2(ii)]{EGAIV},
there exist $\um=\um_i$ and an $\um$-variety $Z_\eta$ which induces $X$
in the sense that $X \xrightarrow{\cong} Z_\eta\otimes_{\um} \ul$. (The latter
is an abbreviation for $Z_\eta\otimes_{\Spec \um} \Spec \ul$.) Moreover,
$X$ is the projective limit of
the $\uk$-schemes $X\otimes_{\um} \um_i$, for the $\um_i$ which contain $\um$.
By \cite[Prop.\,8.13.1]{EGAIV},
after replacing $\um$ with a larger $m_i$
if necessary, there is a $\uk$-morphism $g:Z_\eta\to Y$ which induces the
natural $\uk$-morphism $X \to Y$ in the sense that the latter factors
through $g$. In particular, we obtain a factorization
$X \xrightarrow{\cong} Z_\eta\otimes_{\um} \ul \to Z_\eta \to Y$.

Now we construct $Z$ by approximating $Z_\eta$ with a $\uk$-variety.
Take an integral $\uk$-variety $M_0$ with field of fractions $\um$.
Then $\Spec \um$ is the projective limit of all open subvarieties
$M_i\hookrightarrow M_0$. By
\cite[Thm.\,8.8.2(ii)]{EGAIV}, there exist $i$ and a morphism $Z_i\to M_i$
of finite type such that
$Z_\eta=Z_i\times_{M_i}\Spec\um$; then $Z_\eta$ is the projective limit of
the schemes
$Z_j=Z_i\times_{M_i}M_j$ for $M_j\hookrightarrow M_i$. Obviously, each
morphism $Z_\eta\to Z_i$ is a generic-fibre
embedding. By \cite[Prop.\,8.13.1]{EGAIV}, the $\uk$-morphism $Z_\eta\to Y$
is induced by a morphism $Z\to Y$
for an appropriate choice of $Z=Z_j$, i.e. $Z_\eta\to Y$ factors through
a morphism of $\uk$-varieties $Z\to Y$.
This proves (1).

Now we prove (2). Note that $f_{\ul}$ is the base change obtained from $h: \Spec \ul \to \Spec \um$; hence
$f_{\ul}$ is faithfully flat and $f_{\ul}$ is regular if and only if $h$ is regular. The latter condition is
automatic in characteristic zero. Note also that $f_{\um}$ is regular because it is a \emph{pro-open immersion}
in the sense of \cite[\S2.1]{Tem}
(i.e. $f_{\um}$ is a projective limit of open immersions; in particular, it is injective and
$\cO_Z|_{Z_\eta}=\cO_{Z_\eta}$)
, and $f_{\uk}$ is of finite type, hence it is regular if and only if it is smooth. Since $f=f_{\uk}\circ
f_{\um}\circ f_{\ul}$ and regularity is preserved by composition \cite[Lemma 33.B]{Mat}, we see that $f$ is
regular provided that $f_{\uk}$ is smooth, and clearly then provided that $f_{\uk}$ is smooth on a neighborhood
of $Z_\eta\hookrightarrow Z$. Conversely, suppose that $f$ is regular. By \cite[Lemma 33.B]{Mat}, since
$f_{\ul}$ is faithfully flat, the morphism $Z_\eta\to Y$ is regular. Let $T\subset Z$ be the non-smooth locus of
$f_{\uk}$; then a point $z\in Z$ lies in $T$ if and only if the morphism $f_z:\Spec \cO_{Z,z} \to \Spec
\cO_{Y,f_{\uk}(z)}$ is not regular. Consider $z\in Z_\eta$. Then the local ring of $z$ in $Z_\eta$ is the same
as its local ring in $Z$, because $\cO_{Z_\eta}=\cO_Z|_{Z_\eta}$. Therefore, $f_z$ is regular because $Z_\eta\to
Y$ is regular. So $z\notin T$; thus $Z_\eta$ is disjoint from the closed set $T$ and $f_{\uk}$ is smooth on
$U:=Z\setminus T$, as required.
\end{proof}

\subsection{Application to functorial resolution of singularities}

\begin{corollary}
A desingularization algorithm for algebraic varieties in characteristic zero
is functorial with respect to regular morphisms if and only if it is functorial
with respect to smooth morphisms, ground-field extensions and generic-fibre
embeddings.
\end{corollary}

\begin{remarks}
(1) A \emph{functorial} desingularization algorithm associates to a variety $X$
a sequence of blowings-up $\cF(X)$ with the property that, if $f: X \to Y$ is an
allowed morphism (e.g., a regular morphism in Corollary 4.6), then the
desingularization sequence $\cF(Y)$ pulls back to $\cF(X)$, after perhaps
deleting isomorphisms in the pulled-back sequence \emph{when
$f$ is not surjective}. (For example, if $f$ is an open immersion, then the
centre of a given blowing-up in the resolution sequence for $Y$ may have no
points over $X$. See also \cite[\S\S2.3.3--2.3.6]{Tem2}.)
\smallskip

(2) Suppose that $Y$ is a $\uk$-variety, where $\car \uk = 0$, and that $f: X \to Y$
is a regular morphism of varieties. Consider a finite covering $\{X_i\}$ of $X$
by connected open affine subvarieties. For each $i$, let
$\g_i: X_i \hookrightarrow X$ denote the inclusion and set $f_i := f\circ \g_i:
X_i \to Y$. For each $i$, there is a morphism of finite type $X_i \to \Spec \ul_i$
(by Corollary \ref{concor}) and $f_i$ factors according to Theorem 1.2;
let us denote the
factorization as
\begin{equation*}
X_i \,\cong\, (Z_i)_{\eta_i} \times_{\Spec \um_i} \Spec \ul_i \,
\xrightarrow{(f_i)_{\ul_i}}\, (Z_i)_{\eta_i}
\,\xrightarrow{(f_i)_{\um_i}}\, Z_i \,\xrightarrow{(f_i)_{\uk}}\, Y.
\end{equation*}
There are induced morphisms
\begin{equation}
\coprod_i X_i \, \to \, \coprod_i (Z_i)_{\eta_i} \, \to \,
\coprod_i Z_i \, \to \, Y,
\end{equation}
where $\coprod$ denotes disjoint union. The proof below involves
functoriality with respect to these morphisms. But $\coprod_i (Z_i)_{\eta_i}$ is
not necessarily a variety since it does not necessarily admit a morphism of finite
type to $\Spec$ of a fixed field. It is therefore convenient to extend our
desingularization theorems to schemes that are disjoint unions of varieties.
(We will extend the use of \emph{variety} to include such schemes.)
The desingularization theorems of \cite{BMfunct} extend trivially to this larger
category.
\end{remarks}

\begin{proof}[Proof of Corollary 4.6]
Assume we have a desingularization algorithm that is functorial with respect to smooth morphisms, ground-field
extensions and generic-fibre embeddings. Let $Y$ be a $\uk$-variety, where $\car \uk = 0$, and let $f: X \to Y$
denote a regular morphism of varieties. We have to show that the desingularization sequence for $Y$ pulls back
to that for $X$ (modulo trivial blowings-up in the pull-back sequence; cf.
Remarks 4.7(1)). 
We use the notation of Remarks 4.7(2). 
The morphism $\g: \coprod X_i \to X$ induced by the inclusions $\g_i: X_i
\hookrightarrow X$ is \'etale and surjective. 
By functoriality with respect to smooth (and hence, in particular, \'etale)
morphisms, it is therefore enough to show that the desingularization algorithm commutes with pullback by the
composite of the three morphisms in (5.2). This is true by the assumption.
\end{proof}

\begin{remark}
Another approach to functoriality (and, in particular,
to the assertion of Corollary 4.6) of origin in \cite{BMinv} involves proving a stronger desingularization
theorem where the centres of blowings-up of a variety $X$ are given by the maximum loci of an
upper-semicontinuous \emph{desingularization invariant} (see \cite[\S7]{BMfunct}). Functoriality of the
algorithm with respect to smooth morphisms, ground-field extensions and generic-fibre embeddings then implies
functoriality with respect to regular morphisms, directly by Theorem 1.2.
\end{remark}

Both Corollary 4.6 and Remark 4.8 have analogues for desingularization
of marked ideals that can be obtained in the same way.

\section{Functoriality of desingularization of a marked ideal
with respect to generic-fibre embeddings}
In this section we prove that the desingularization algorithm for marked
ideals \cite[\S5]{BMfunct} is functorial with respect to generic-fibre
embeddings (Proposition 5.3 below). Theorem 1.3 then follows from
Corollary 4.6 together with functoriality with respect to ground-field extensions
and with respect to smooth morphisms (\cite[\S7]{BMfunct}; see Remark 5.1 below).

Let $\pi: Z \to T$ denote a dominant morphism of $\uk$-varieties, where
$\uk$ is a field of characteristic zero. Let $\psi: Z_\eta \to Z$ denote
the embedding of the generic fibre of $\pi$, as in Section 2. ($Z_\eta$
is an $\um$-variety, where $\um = K(T)$.) If $\cI \subset \cO_Z$ is an
ideal (i.e., a coherent sheaf of ideals), let $\cI_\eta \subset \cO_{Z_\eta}$
denote the inverse image (pullback) $\psi^*(\cI)$. Then $\cI_\eta$ is a
coherent sheaf of ideals on $Z_\eta$.

Every subvariety of $Z_\eta$ is of the form $C_\eta$, where $C$ is a subvariety
of $Z$ such that $\pi|_C: C \to T$ is dominant. Moreover $C_\eta$ is smooth
if and only if there is an open subset of $T$ over which $C$ (and also $\pi|_C$)
is smooth (Section 2). Let $\Bl_C Z \to Z$ denote the blowing-up with centre a
subvariety $C$ of $Z$. By functoriality of blowing-up with respect to flat
base extension,
\begin{equation}
(\Bl_C Z)_\eta\,=\,\Bl_{C_\eta}Z_\eta
\end{equation}
(where we understand $C_\eta = \emptyset$ and $\Bl_{C_\eta}Z_\eta = Z_\eta$
if $\pi|_C$ is not dominant).

\subsection{Marked ideals}
A \emph{marked ideal} $\ucI$ is a quintuple $\ucI = (Z,N,E,\cI,d)$, where:
$Z \supset N$ are smooth varieties,
$E = \sum_{i=1}^s H_i$ is a simple normal crossings divisor on $Z$ which is
tranverse to $N$ and \emph{ordered} (the $H_i$ are smooth hypersurfaces
in $Z$, not necessarily irreducible, with ordered index set as indicated),
$\cI \subset \cO_N$ is an ideal, and $d \in \IN$.
The \emph{cosupport} of $\ucI$,
$$
\cosupp \ucI := \{x \in N:\, \ord_x\cI \geq d\}\,.
$$
We say that $\ucI$ is of \emph{maximal
order} if $d = \max\{\ord_x \cI: x \in \cosupp \ucI\}$. The \emph{dimension}
$\dim \ucI$ denotes $\dim N$.

A blowing-up $\s: Z'= \Bl_C Z \to Z$ (with smooth centre $C$) is
\emph{admissible} for $\ucI$ if
$C \subset \cosupp \ucI$, and
$C$, $E$ have only normal crossings.
The \emph{(controlled) transform} of $\ucI$ by an admissible blowing-up
$\s: Z' \to Z$ is the marked ideal $\ucI' = (Z',N',E',\cI',d'=d)$,
where
$N'$ is the strict transform of $N$ by $\s$,
$E' = \sum_{i=1}^{s+1} H_i'$, where $H_i'$ denotes the strict transform
of $H_i$, for each $i=1,\ldots,s$, and $H'_{s+1} := \s^{-1}(C)$ (the
exceptional divisor of $\s$, introduced as the last member of $E'$), and
$\cI' := \cI_{\s^{-1}(C)}^{-d}\cdot \s^*(\cI)$ (where $\cI_{\s^{-1}(C)}
\subset \cO_{N'}$ denotes the ideal of $\s^{-1}(C)$).
In this definition, note that $\s^*(\cI)$ 
is divisible by $\cI_{\s^{-1}(C)}^d$ and $E'$ is a normal crossings divisor transverse to $N'$, because $\s$ is
admissible.

We define a \emph{resolution of singularities} of a marked ideal $\ucI$ as
a finite sequence of admissible blowings-up after which $\cosupp \ucI = \emptyset$.

Let $\ucI$ be a marked ideal as above, and let $\vp: Y \to Z$ denote a regular
morphism. We define the \emph{inverse image} (or \emph{pullback}) $\vp^*(\ucI)$
as the marked ideal $(Y, \vp^{-1}(N), \vp^{-1}(E), \vp^*(\cI), d)$ (where
$\vp^{-1}(E)$ inherits the ordering of $E$). If $\psi: Z_\eta \to Z$ is a
generic-fibre embedding as above, then $\vp^*(\ucI) = \ucI_\eta$, where
the latter denotes the marked ideal $(Z_\eta, N_\eta, E_\eta, \cI_\eta, d)$.
($N_\eta$ is empty (so $\cosupp \ucI_\eta = \emptyset$) unless $\pi|_N$ is
dominant, and $E_\eta$ is empty if no component of $E$ dominates $T$.)

\begin{remark}
There are two proofs of functoriality of the desingularization
algorithm for a marked ideal with respect to \'etale or smooth
morphisms in \cite[\S7]{BMfunct}. The first proof comes from Kollar
\cite[Prop.\,3.37]{Ko}. For this argument, we assume that our marked ideals are
equidimensional and that the smooth or regular morphisms considered
have constant relative
dimension; it seems inconvenient to carry out the proof without these
assumptions. The proof is by induction on dimension. It uses functoriality in the
inductive step in a way that necessitates working with marked ideals
$\ucI = (Z,N,E,\cI,d)$ where $Z$ may have several components (cf. Remarks 4.7(2)).
Although the
blowings-up of different components are independent, a functorial algorithm
depends on which we take first, second, etc. (In dimension $1$, the algorithm
dictates blowing up the points of maximum order of $\cI$ at each step.)

The second proof of functoriality with respect to smooth morphisms
involves proving the stronger desingularization theorem
\cite[Thm. 7.1]{BMfunct}, where the centres of blowing up are given by the
maximum loci of an upper-semicontinuous \emph{desingularization invariant}
(cf. Remark 4.8).
The invariant is defined by induction on $\dim \ucI$. This functorial
desingularization algorithm does not require an equidimensionality
assumption on the marked ideals, and applies to smooth or regular morphisms
that are not necessarily of constant relative dimension.
\end{remark}

\subsection{Test transformations and equivalence}
Let $\ucI$ denote a marked ideal as above. Sequences of test transformations are
introduced to test for invariance of local numerical characters of $\ucI$ (see
\cite[\S6]{BMfunct}).
\emph{Test transformations} are transformations of a marked ideal
by morphisms of three possible
kinds: admissible blowings-up, projections from products with an affine line, and
exceptional blowings-up:

\emph{Product with a line}.
Let $Z' := Z \times \IA^1$, and let
$\pi:\, Z' \to Z$ denote the projection. We define the \emph{transform}
$\ucI'$ of $\ucI$ by $\pi$ as the marked ideal $\ucI' = (Z',N',E',\cI',d'=d)$,
where $N' := \pi^{-1}(N)$, $\cI' := \pi^*(\cI)$, but $E' =
\sum_{i=1}^{s+1} H_i'$, where $H_i' := \pi^{-1}(H_i)$, for each $i = 1,\ldots,s$,
and $H'_{s+1}$ denotes the \emph{horizontal divisor} $D := Z \times \{0\}$
(included as the last member of $E'$).

\emph{Exceptional blowing-up}.
A blowing-up $\s: Z' \to Z$
is called an \emph{exceptional blowing-up} for $\ucI$ if its centre
$C$ is an intersection $H_i \cap H_j$ of distinct hypersurfaces $H_i, H_j
\in E$.
We define the \emph{transform} $\ucI' = (Z',N',E',\cI',d')$ of $\ucI$ by
$\s$ as the marked ideal in the same way as for an admissible blowing-up.
(In the case of an exceptional blowing-up, $N'=\s^{-1}(N)$ and $\cI' = \s^*(\cI)$.)

A \emph{test sequence} for $\ucI_0 = \ucI$ means a sequence of morphisms
\begin{equation*}
Z = Z_0 \stackrel{\s_1}{\longleftarrow} Z_1 \longleftarrow \cdots
\stackrel{\s_{t}}{\longleftarrow} Z_{t}\,,
\end{equation*}
where each successive $\s_{j+1}$ is either an admissible blowing-up,
the projection from a product with a line, or an exceptional blowing-up.

We say that two marked ideals $\ucI$ and $\ucI_1$ (with the same ambient
variety $Z$ and the same normal crossings divisor $E$) are \emph{equivalent}
if they have the same test sequences (i.e., every test sequence for one is a
test sequence for the other).

\begin{lemma}
Let $Z_\eta \to Z$ denote a generic-fibre embedding as above. Suppose that
$\ucI = (Z,N,E,\cI,d)$ and $\ucJ = (Z,P,E,\cJ,e)$ are marked ideals on $Z$.
If $\ucI$ is equivalent to $\ucJ$, then $\ucI_\eta$ is equivalent to $\ucJ_\eta$.
\end{lemma}

\begin{proof}
This follows directly from the definitions, together with the fact that
any test sequence for $\ucI_\eta$ lifts to a test sequence for $\ucI$ over
some neighbourhood of $Z_\eta$ in $Z$ (cf. (5.1)) and the
fact that if $b \in N_\eta$ and $a = \psi(b)$, then $\ord_a \cI = \ord_b \cI_\eta$
(see Section 2).
\end{proof}

\subsection{Functoriality with respect to generic-fibre embeddings}

\begin{proposition}
Let $\psi: Z_\eta \to Z$ denote a generic-fibre embedding and let $\ucI$
denote a marked ideal on $Z$, as above. Then the sequence of blowings-up
involved in the desingularization algorithm \cite[\S5]{BMfunct} for $\ucI_\eta$
is the pullback of the desingularization sequence for $\ucI$.
\end{proposition}

\begin{remark}
In the blowing-up sequence for $\ucI$, any centre of blowing up that does
not dominate $T$ pulls back to an empty centre, so that the corresponding
blowing-up over $Z_\eta$ is the identity morphism.
\end{remark}

\begin{proof}[Proof of Proposition 5.3]
The proof consists simply of following that in \cite[\S5]{BMfunct}
step-by-step. The proof is by induction on $\dim \ucI$.
We will not go through the entire process. The proof (or the algorithm)
as presented in \cite[\S5]{BMfunct} has two main steps, each of which
involves an important construction: in Step I, passage from a marked ideal
$\ucI = (Z,N,E,\cI,d)$ of maximal order to a
\emph{coefficient ideal} $\ucC(\ucI)$ on an open subset $U$ of $Z$, to
decrease the dimension by 1 (for induction), and in Step II, passage from
a general marked ideal $\ucI$ to a \emph{companion ideal} $\ucG(\ucI)$ of
maximal order, to reduce to Step I.

It is easy to see that $\ucG$ commutes with
$\ucI \mapsto \ucI_\eta$ (by the definition of $\ucG$;
we leave the details to the reader). For $\ucC$,
commutativity with respect to $\ucI \mapsto \ucI_\eta$ is true only on the
level of equivalence classes. This is proved in the following subsection.
Our proposition follows from these two results.
\end{proof}

\begin{remark}
Step II in \cite[\S5]{BMfunct} involves proving that the equivalence class
of $\ucG(\ucI)$ depends only on the equivalence class (and dimension) of $\ucI$.
This result is proved using the fact that two local numerical characters of
a marked ideal, $\ord_a \cI/d$ and $\ord_{H,a} \cI/d$, $H \in E$
(where $\ord_{H,a}$ denotes the order along $H$) are invariants of the
equivalence class. In \cite{Hid}, Hironaka proposes to prove Theorem 1.3 above
using a weaker notion of equivalence where test sequences involve only
admissible blowings-up and product with an affine line. Although $\ord_a \cI/d$
is an invariant of the weaker equivalence class, $\ord_{H,a} \cI/d$ is not
\cite[Ex.\,5.14]{BMda1}.
\end{remark}

\subsection{Coefficient ideals}
Let $\ucI = (Z,N,E,\cI,d)$ denote a marked ideal as above. Recall from Section 1
that the derivative ideal $\cD(\cI)$ is the image of the natural morphism
$\Der_N \times \cI \to \cO_N$. Let $\cD_E(\cI) \subset \cO_N$ denote the ideal
generated by all local sections of $\cI$ and all derivations that preserve the
ideal $\cI_E$ of $E$. Higher-derivative ideals are defined inductively by
$$
\cD^{j+1}_E(\cI) := \cD_E(\cD^j_E(\cI)), \quad j= 1,\ldots.
$$
We define marked ideals
$$
\ucD_E^j(\ucI) := (M,N,E,\cD^j_E(\cI),d-j),\quad j = 1,\ldots,d-1,
$$
and
$$
\ucC_E^k(\ucI) := \sum_{j=0}^k \ucD_E^j(\ucI), \quad k \leq d-1;
$$
$\ucC_E^k(\ucI)$ is a marked ideal
$\left(M,N,E, \cC_E^k(\ucI), d_{\ucC_E^k(\ucI)}\right)$. ($\ucC_E^k(\ucI)$
is a weighted sum of marked ideals; see \cite[\S3.3]{BMfunct}.)
The marked ideals $\ucI$ and $\ucC_E^k(\ucI)$, $k \leq d-1$, are equivalent
\cite[Cor.\,3.1]{BMfunct}.

Suppose that $\ucI = (Z,N,E,\cI,d)$ is of maximal order. Then every point
of $\cosupp \ucI$ has an open neighbourhood $U \subset Z$ in which $\ucI$
has a \emph{maximal
contact hypersurface} $P \subset N|_U$ \cite[\S4]{BMfunct}. The corresponding
\emph{coefficient ideal} is defined as
$$
\ucC_{E,P}(\ucI) := \left(U,P,E, \cC_E^{d-1}(\ucI)|_P,
d_{\ucC_E^{d-1}(\ucI)}\right).
$$
It follows that $\ucC_{E,P}(\ucI)$ is equivalent to $\ucI|_U$
\cite[Cor.\,4.1]{BMfunct}.

\begin{lemma}
Suppose that $\psi: Z_\eta \to Z$ is a generic fibre embedding as above.
Then $\ucC_{E,P}(\ucI)_\eta = \psi^*\ucC_{E,P}(\ucI)$ is equivalent to
$\ucC_{E_\eta,P_\eta}(\ucI_\eta)$.
\end{lemma}

\begin{proof}
This is immediate from the preceding equivalence and Lemma 5.2.
\end{proof}

\begin{remark}
Lemma 5.6 
can be understood also in another way. The marked ideal $\ucC_{E,P}(\ucI)$ is equivalent to a smaller
coefficient ideal $\ucC_{z,P}(\ucI)$ defined using only derivatives in the normal direction to the hypersurface
$P \subset N$ (i.e., derivatives with respect to a local generator $z$ of the ideal of $P$ in $\cO_N$)
\cite[Ex.\,4.4(1)]{BMfunct}. For this variant of the coefficient ideal, $\psi^*\ucC_{z,P}(\ucI) =
\ucC_{z|_{N_\eta},P_\eta}(\ucI_\eta)$. (Derivatives along elements of $\um$ that are transcendental over $\uk$
are not explicitly involved; cf. Remark 1.4.) Lemma 5.6 
follows also from the latter, \cite[Ex.\,4.4(1)]{BMfunct}, and Lemma 5.2.
\end{remark}

\section{Functoriality of desingularization of a variety}
We begin with a precise statement of Theorem 1.1.


\begin{theorem}
Given a variety $X$ over a field $\uk$ of characteristic zero,
there is finite sequence of blowings-up $\s_{j+1}: X_{j+1} \to X_j$
with smooth centres,
\begin{equation}
X = X_0 \stackrel{\s_1}{\longleftarrow} X_1 \longleftarrow \cdots
\stackrel{\s_{t}}{\longleftarrow} X_{t}\,,
\end{equation}
such that:
\begin{enumerate}
\item
$X_t$ is smooth and the exceptional divisor in $X_t$ has only normal crossings.
\item
All centres of blowing up are disjoint from the preimages of $X\setminus \Sing X$.
\item
The resolution morphism $\s_X: X_t \to X$ given by the composite of the
$\s_j$ (or the entire sequence of blowings-up (6.1)) is associated
to $X$ in a way that is functorial with respect to regular morphisms.
(See Remarks 4.7(1).)
\end{enumerate}
\end{theorem}

This theorem can be proved with the following stronger version of the condition (2): For each $j$, let $C_j
\subset X_j$ denote the centre of the blowing-up $\s_{j+1}: X_{j+1} \to X_j$. Then either $C_j \subset \Sing
X_j$ or $X_j$ is smooth and $C_j$ lies in the support of the exceptional divisor of $\s_1\circ\cdots\circ\s_j$.
In fact, we prove Theorem 6.1 together with the following addendum (where
(3), (4) should again be understood modulo trivial blowings-up as in 
Remarks 4.7(1)).

\begin{addendum}
Given any embedding (i.e., closed immersion) $X|_U \hookrightarrow Z$, 
where $U$ is an open subset
of $X$ and $Z$ is smooth, there is a sequence of blowings-up
$\tau_{j+1}: Z_{j+1} \to Z_j$,
\begin{equation}
Z = Z_0 \stackrel{\tau_1}{\longleftarrow} Z_1 \longleftarrow \cdots
\stackrel{\tau_{t}}{\longleftarrow} Z_{t}\,,
\end{equation}
which satisfies the following conditions. Set $Y_0 := X|_U$.
For each $j$, let $C_j$ denote
the centre of $\tau_{j+1}$, let $E_{j+1}$ denote the exceptional divisor
of $\tau_1\circ\cdots\circ\tau_{j+1}$, and define $Y_{j+1}$
inductively as the strict transform of $Y_j$ by $\tau_{j+1}$. Then:
\begin{enumerate}
\item
Each $C_j$ is smooth and has only normal crossings with respect to $E_j$.
\item
For each $j$, either $C_j \subset \Sing Y_j$ or $Y_j$ is smooth and
$C_j \subset Y_j \cap \supp E_j$.
\item
Each $X_j|_U = Y_j$ and, over $U$, the resolution sequence (6.1) is given
by the restriction of (6.2) to the $Y_j$.
\item The sequence of blowings-up (6.2) is associated to $X|_U \hookrightarrow Z$
in a way that is functorial with respect to regular morphisms (of the ambient
variety $Z$).
\end{enumerate}
\end{addendum}

A weaker version of Theorem 6.1 (and the Addendum) can be obtained directly from Theorem 1.3 applied to the
marked ideal $\ucI = (Z,Z,\emptyset,\cI_X,1)$, where $X \hookrightarrow Z$ is a (local) embedding of $X$ in a
smooth $\uk$-variety $Z$ and $\cI_X \subset \cO_Z$ is the ideal of $X$ (see \cite[\S1.1]{BMfunct}. In this
version, the resolution sequence (6.1) is given by restricting the blowing-up sequence provided by Theorem 1.3
to the successive strict transforms $X_j$ of $X$ --- the intersections of the centres $C_j$ with the $X_j$ will
not necessarily be smooth, nor will condition (2) of the Addendum necessarily hold.

Theorem 6.1 (including the Addendum) as stated except for weaker functoriality
conditions --- with respect to smooth morphisms and base-field extensions ---
is proved in \cite{BMinv, BMfunct} (under the tacit assumption that the smooth
morphisms have constant relative dimension. We show how to remove this restriction
in \S6.3 below.) The proof involves the \emph{Hilbert-Samuel
function} $H_{X,a} \in \IN^{\IN}$, $a \in X$, and desingularization of an associated
marked ideal that we call a \emph{presentation} of the Hilbert-Samuel function.
We will use commutativity of a presentation with respect to generic-fibre embeddings
(Proposition 6.3 below) together with Theorem 1.2 and Remark 4.8
to deduce Theorem and Addendum 6.1 in full (see \S6.4).

\subsection{The Hilbert-Samuel function}
If $a$ is a closed point, then the \emph{Hilbert-Samuel function} $H_{X,a}$
is defined as
$$
H_{X,a}(k) := \length \frac{\cO_{X,a}}{\fm_{X,a}^{k+1}}\,,\quad k \in \IN\,.
$$
Thus $H_{X,a} \in \IN^{\IN}$.
We can extend the definition to arbitrary points of $X$ so that $a \mapsto
H_{X,a}$ will be upper-semicontinuous (where the set of sequences $\IN^{\IN}$ is
totally-ordered lexicographically):

In general, define $\La: \IN^{\IN} \to \IN^{\IN}$ by
$$
\La(F)(k) = \sum_{j=0}^k F(j)\,,\quad \quad k \in \IN\,,
$$
where $F \in \IN^{\IN}$. Define $\La^j: \IN^{\IN} \to \IN^{\IN}$,
$j \geq 1$, inductively by $\La^j(F)
= \La(\La^{j-1}(F))$. Suppose that $R$ is a Noetherian local ring with maximal ideal
$\fm$. Define $H^{(0)}(R) \in \IN^{\IN}$ by
$$
H^{(0)}(R)(k) := \length \frac{R}{\fm^{k+1}}\,,\quad k \in \IN\,.
$$
For each $j \in \IN$, let $H^{(j)}(R)$ denote $\La^j(H^{(0)}(R))$. If $a \in X$, then
we define $H_{X,a}^{(j)} := H^{(j)}(\cO_{X,a})$, for all $j \in \IN$, and
we let $H_{X,a}$ denote $H_{X,a}^{(l)}$, where $l$
denotes the transcendence degree of the residue field $\ka(a)$ over $\uk$
(i.e., the dimension of the closure of $a$).

The Hilbert-Samuel function $H_{X,a}$ determines the \emph{minimal embedding
dimension} $e_{X,a}$ of $X$ at $a$ (in a smooth affine $\uk$-variety):
$e_{X,a} = H_{X,a}(1) - 1$.

The Hilbert-Samuel function $H_{X,\cdot}: X \to \IN^{\IN}$
has the following basic properties, established
by Bennett \cite{Be} (see \cite{BMjams, BMinv} for simple proofs):
(1) $H_{X,\cdot}$ distinguishes smooth and singular points. (2) $H_{X,\cdot}$ is
(Zariski) upper-semicontinuous. (3) $H_{X,\cdot}$ is \emph{infinitesimally
upper-semicontinuous} (i.e., $H_{X,\cdot}$ cannot increase after blowing-up
with centre on which it is constant). (4) Any decreasing sequence in the value
set of the Hilbert-Samuel function stabilizes.

\subsection{Presentation of the Hilbert-Samuel function}
The Hilbert-Samuel function $H_{X,a}$ is a local invariant that
plays the same role with respect to strict transform of a variety $X$
as the order plays with respect to (weak) transform of a (marked) ideal.
More precisely, for all $a \in X$, there is an embedding $X|_U \hookrightarrow
Z$, where $U$ is a neighbourhood of $a$ and $Z$ is smooth,
and a marked ideal $\ucI = (Z,N,\emptyset,\cI,d)$ which has the same
\emph{test sequences} as $\uX := (Z,\emptyset,X|_U,H)$, where $H = H_{X,a}$. (We
define a \emph{test sequence} for $\uX = (M,E,X,H)$
by analogy with that for a marked
ideal (\S5.2), but where a blowing-up $\s: Z' \to Z$ with smooth centre $C$
is \emph{admissible} if $C \subset \supp \uX := \{x\in X: H_{X,x} \geq H\}$ and
$C,\,E$ have only normal crossings, and where $X$ transforms by \emph{strict
transform}.) We call $\ucI$ a \emph{presentation} of $H_{X,\cdot}$ at $a$.

A construction of a presentation is given in \cite[Ch.\,III]{BMinv} and some
familiarity with the latter will be needed to understand the results of this
section in detail. The essential point needed is that $\cI$ is generated by
suitable powers of a special system of generators of $\ucI_X \subset \cO_Z$
at $a$ which
is determined by the \emph{vertices} of the \emph{diagram of initial exponents}
$\fN(\cI_{X,a})$ with respect to a local coordinate system for $Z$ at $a$
(see \S6.5 below).

We can choose a presentation $\ucI = (Z,N,\emptyset,\cI,d)$ of $H_{X,\cdot}$
at $a$ so that $Z$ is a smooth minimal embedding variety for $X$ at $a$. Given
$\ucI$, there is an \'etale morphism $\vp: Z'\to Z$ onto a neighbourhood of $a$
such that $\vp^*(\ucI)$ is equivalent to a marked ideal $\ucJ = (Z',N',\emptyset,
\cJ,e)$ of maximal order and codimension zero (i.e., $N'=Z'$)

\subsection{Functoriality with respect to smooth morphisms}
If $a$ is a maximum point of $H_{X,\cdot}$ and $\ucI$ is a presentation of
$H_{X,\cdot}$ at $a$, then the corresponding maximal value of
$H_{X,\cdot}$ decreases after desingularization of $\ucI$. This is the main
point of a presentation of the Hilbert-Samuel function, needed to prove the
strong desingularization theorem for a variety using functorial desingularization
of a marked ideal together with the basic properties of the Hilbert-Samuel
function in \S6.1 above.

In \cite[\S7]{BMfunct}, versions of Theorem 6.1, where the functoriality
assertion is with respect to smooth morphisms and base-field extensions,
are proved using a presentation of the Hilbert-Samuel function, again
following either of the two schemes recalled in Remark 5.1 above.

However, there are equidimensionality issues for smooth morphisms that
are not treated in previous works, even using a desingularization invariant.
(These issues were raised in a letter from Ofer Gabber to the third author.)
We can deal with them only using the second of the two methods in
Remark 5.1, which again involves proving a stronger
desingularization theorem where the centres of blowing up are given by the
maximum loci of a desingularization invariant $\inv_X$ defined inductively over
a sequence of admissible blowings-up. Since the marked ideal $\ucJ$ above is
of maximal order, the desingularization invariant $\inv_{\ucJ}$ for $\ucJ$ is
a finite sequence whose first term is $1$ throughout the cosupports of the
successive transforms of $\ucJ$ (see \cite[\S7.2]{BMfunct}); $\inv_X$ is
defined at the corresponding points of $X$ and its successive strict transforms
by replacing this first term by $H_{X,a}$.
For details  we refer to \cite{BMinv, BMda1, BMfunct}. (It is important to
begin with a presentation $\ucJ$ as above so that the desingularization
invariant will be independent of the choice of a local embedding variety for $X$.)

As shown in \cite{BMfunct},
blowing up with centre $=$ maximum locus of $\inv_X$ gives an algorithm
for resolution of singularities of arbitrary $X$, functorial with respect
to smooth morphisms of constant relative dimension.

In order to prove functoriality with respect to arbitrary smooth morphisms, we
first note that (in the notation of \S6.1), $\La^k(H_{X,a}) = H_{X\times \IA^k, (a,0)}$.
(See also \S6.5 below.)
Moreover a presentation $\ucI = (Z,N,\emptyset,\cI,d)$ of $H_{X,\cdot}$ at $a$
induces a presentation of
$H_{X\times \IA^k, \cdot}$ at $(a,0)$ by pull-back by the projection $Z\times \IA^k
\to Z$.

Suppose that $X$ is locally equidimensional. Define a modified invariant
$\inv^*_X$ by replacing the first term $H_{X,x}$ of $\inv_X(x)$ at each
(closed) point $x$ by
$\La^{d-q}(H_{X,x})$, where $d = \dim X$ and $q = q(x) := \dim \cO_{X,x}$.
Then blowing up with centre $=$ maximum locus of $\inv^*_X$ gives an algorithm
for resolution of singularities of locally equidimensional varieties $X$,
\emph{functorial with respect to arbitrary smooth morphisms}.

We can use the preceding idea together with a suggestion of Gabber (in the letter
cited above) to define a modified invariant $\inv^\#_X$ such that blowing up with
centre $=$ maximum locus of $\inv^\#_X$ gives an algorithm for \emph{resolution of
singularities of arbitrary varieties $X$, functorial with respect to
arbitrary smooth morphisms}: Let $\#(x)$ denote the number of different dimensions
of %
irreducible components of $X$ at $x$. Let $q(x)$ denote the smallest dimension of the irreducible components.
Define $\inv^\#_X$ by replacing the first term $H_{X,x}$ of $\inv_X(x)$ at each (closed) point $x$ by the pair
$(\#(x),\, \La^{d-q}(H_{X,x}))$, where $d = \dim X$ and $q = q(x)$. It is easy to see that a marked ideal is a
presentation of the Hilbert-Samuel function at $x$ if and only if it is a presentation of 
$(\#(\cdot),\,\La^{d-k}(H_{X,\cdot}))$. The assertion follows.

\begin{remark}
The fact that the invariants $H_{X,\cdot}$ and $(\#(\cdot),\,\La^{d-k}(H_{X,\cdot}))$
share a common presentation at every point means, in particular, that every component
of a constant locus of one of these invariants is also a component of a constant locus
of the other. It follows that the maximum loci of the two invariants $\inv_X$ and
$\inv^\#_X$ are each unions of closed components of constant loci of the invariant
$\inv_X$, but not necessarily of 
the same closed components in each case --- i.e., the
order in which we blow up these components may not be the same. The invariant
$(\#(\cdot),\,\La^{d-k}(H_{X,\cdot}))$ is contrived to force us 
to blow up components
in an order that gives functoriality with respect to arbitrary smooth morphisms.
\end{remark}

\subsection{Functoriality with respect to generic-fibre embeddings}
In order to deduce Theorem 6.1 (and its Addendum) with the full version of
functoriality, i.e., with respect to regular morphisms in general, using Theorem 1.2
and Remark 4.8, it is now again enough
to prove functoriality with respect to generic-fibre embeddings. For the latter,
because of Lemma 5.2 and Proposition 5.3, it is enough to prove Proposition 6.3
following.

Let $\psi: X_\eta \to X$ denote a generic-fibre embedding, corresponding to
a dominant morphism of $\uk$-varieties $\pi: X \to T$, where $T$ is integral.
Let $b$ denote a closed point of $X_\eta$ and let $a = \psi(b) \in X$.
Then there is a neighbourhood $U$ of $a$ in $X$ so that $X|_U$ embeds in a smooth
$\uk$-variety $Z$ such that $\pi$ extends to a (dominant) morphism
$Z \to T$. We can choose $U$ and $Z$ such that $Z$ is a minimal embedding
variety for $X$ at $a$.

For simplicity of notation, we will write simply $X$ instead of $X|_U$, and
$X_\eta$ for the generic fibre of the latter. Then $X_\eta = X \times_Z Z_\eta$
and there is a commutative diagram
$$
\begin{CD}
X_\eta @>\psi >> X \\
@VVV                   @VVV \\
Z_\eta    @>>>  Z\\
\end{CD}
$$
where the horizontal arrows are the generic-fibre embeddings and the right
(respectively, left) vertical arrow is a morphism of $\uk$-varieties
(respectively, $\um$-varieties, where $\um = K(T)$).

\begin{proposition}
With the notation preceding, there is a neighbourhood $V$ of $a$ in $Z$
and a presentation $\ucI = (Z|_V,N,\emptyset,\cI,d)$ of $H_{X,\cdot}$ at $a$
such that $\ucI_\eta$ is a presentation of $H_{X_\eta,\cdot}$ at $b$.
\end{proposition}

More precisely, we claim that, for suitable local coordinates for $Z$, the
presentation constructed in \cite[Ch.\,III]{BMinv} has the required properties.
Since we do not want to repeat the construction, we will give only the new
ingredients needed by a reader who is familiar with the latter to verify our
claim in a straightforward way.

\subsection{The diagram of initial exponents}
The construction of a presentation in \cite[Ch.\,III]{BMinv} depends on the way
that the Hilbert-Samuel function can be computed using the diagram of initial
exponents of the ideal $\cI_X \subset \cO_Z$ of $X$ with respect to local
coordinates for $Z$ at a point of $X$.

Consider a ring of formal power series $R = K\lbr X \rbr = K\lbr X_1,\dots,
X_n\rbr$ over a field $K$. If $\al = (\al_1,\ldots,\al_n) \in \IN^n$, put
$|\al| = \al_1+\ldots +\al_n$. We totally order $\IN^n$ by using the lexicographic
ordering of $(n+1)$-tuples $(|\al|,\al_1,\ldots,\al_n)$. Consider
$F = \sum_{\al \in \IN^n} F_{\al}X^{\al} \in K\lbr X \rbr$, where
$X^{\al}:=X_1^{\al_1}\cdots X_n^{\al_n}$. Let $\supp F := \{\al: F_{\al}\neq 0\}$.
The \emph{initial exponent} $\exp F$ means the smallest element of $\supp F$.
($\exp F := \infty$ if $F=0$.)

Let $I$ be an ideal in $R$. The \emph{diagram of initial exponents}
$\fN(I) \subset \IN^n$ is defined as
$$
\fN(I) := \{\exp F: F \in I\setminus\{0\}\}.
$$
Clearly, $\fN(I) + \IN^n = \fN(I)$. It follows that there is a smallest
finite subset $\fV$ of $\fN(I)$ (the \emph{vertices} of $\fN(I)$) such that
$\fN(I) = \fV + \IN^n$. ($\fV = \{\al \in \fN(I): (\fN(I)\setminus\{\al\}) + \IN^n
\neq \fN(I)\}$.)

Given $\fN \subset \IN^n$ such that $\fN + \IN^n = \fN$, let $H_{\fN} \in
\IN^{\IN}$ denote the function
$$
H_{\fN}(k) = \#\{\al \in \IN^n\setminus\fN: |\al| \leq k\}\,,\quad k \in \IN
$$
(where $\#S$ denotes the number of elements in a finite set $S$).
Then $H^{(0)}(R/I) = H_{\fN(I)}$ (see \S6.1). It is easy to see that, if
$\fN$ is a product $\fN = \IN^p \times \fN^*$, then $H_{\fN} =
\La^p(H_{\fN^*})$.

Suppose that $(x_1,\ldots,x_n)$ is a coordinate system (system of parameters) on
an open subset $U$ of $Z$. If $c$ is a closed point in $U$, then there is a
unique isomorphism $\wcO_{Z,c} \xrightarrow{\cong} \ka(c)\lbr X_1,\ldots,X_n\rbr$
such that each $x_i \mapsto x_i(c) + X_i$, where $x_i(c)$ denotes the image of
$x_i$ in the residue field $\ka(c) = \cO_{Z,c}/\fm_{Z,c}$. If $\cI_X \subset \cO_Z$
denotes the ideal of $X$, then the \emph{diagram of initial exponents}
$\fN(\cI_{X,c})$ of $\cI_{X,c}$ \emph{with respect to the given coordinate
system} denotes $\fN(I)$, where $I \subset \ka(c)\lbr X\rbr$ is the ideal induced
by $\cI_{X,c}$.

We totally order $\{\fN \subset \IN^n: \fN + \IN^n = \fN\}$ by giving each $\fN$ the
lexicographic order of the sequence of its vertices (in increasing order). Each point
of $X$ admits a coordinate neighbourhood in $Z$ in which the associated diagram
$\fN(\cI_{X,c})$ can be extended to arbitrary points so that $c \mapsto \fN(\cI_{X,c})$
is upper-semicontinuous.

\subsection{Proof of Proposition 6.3}
By restricting $Z$ and $T$ to suitable affine open neighbourhoods of the point $a$
and its image in $T$, we can assume (by the Jacobian criterion for smoothness) that:
\begin{enumerate}
\item
$T$ is a subvariety $V(P)$ of $\IA_{\uk}^{p+q}$ determined by an ideal $(P) \subset
\uk[y,z] = \uk[y_1,\ldots,y_p,z_1,\ldots,z_q]$ generated by polynomials $P_1(y,z),
\ldots,P_q(y,z)$, where the determinant $J_P$ of the Jacobian matrix
$\p P/\p z = \left(\p P_i/\p z_j\right)$ is nonvanishing on $T$.
\item
$Z = V(P,G) \subset \IA_{\uk}^{n+m}$, where $n\geq p,\,m\geq q$, $(P,G)$ is an
ideal in $\uk[x,w]$,
\begin{align*}
x = (u,v) &= (u_1,\dots, u_p, v_1,\ldots, v_{n-p}),\\
w = (s,t) &= (s_1,\dots, s_q, t_1,\ldots, t_{m-q}),
\end{align*}
$(P,G)$ is generated by $P_1(u,s),\ldots,P_q(u,s)$ (from (1)) together with
polynomials $G_1(x,w),\ldots,G_{m-q}(x,w)$, and the determinant $J_{(P,G)}$ of
$\p (P,G)/\p (s,t)$ is nonvanishing on $Z$.
\item
The morphism $Z \to T$ is induced by the inclusion $\uk[y,z] \hookrightarrow
\uk[x,w]$ given by $u=y,\,s=z$.
\end{enumerate}

It follows that
\begin{enumerate}
\item[(4)]
$Z_\eta = V(G_\eta) \subset \IA_{\um}^{(n-p)+(m-q)}$, where
$(G_\eta) \subset \um[v,t] = K(T)[v,t]$ is the ideal generated by the polynomials
$G_{j,\eta}(v,t)$  which are induced by the $G_j(u,v,s,t)$, $j=1,\ldots,m-q$.
\end{enumerate}
Since $J_{(P,G)} = J_P\cdot J_G$, where $J_G = \det\left(\p G/\p t\right)$, we see that
$\det\left(\p G_\eta/\p t\right)$ is nonvanishing on $Z_\eta$.

Therefore, $y=(y_1,\ldots,y_p)$, $x = (x_1,\ldots,x_n) = (u_1,\ldots,u_p,v_1,
\ldots,v_{n-p})$ and $v = (v_1,\ldots,v_{n-p})$ (respectively) induce local coordinates
(regular parameters) on $T$, $Z$ and $Z_\eta$ (respectively).

Let $W$ denote the closure of $a = \psi(b)$ in $X$. Then there is an open subset $V$
of $W$ on which the projection to $T$ is \'etale (see \S2), so that $y=(y_1,\ldots,y_p)$
is a system of coordinates on $V$. Given a closed point $c$ of $V$, let $X_{\pi(c)}$
denote the fibre $X \times_T \pi(c)$ over $\pi(c)$ and let $\cI_{X_{\pi(c)}}$ denote the
ideal of $X_{\pi(c)} \subset Z_{\pi(c)}$. Let $\fN(\cI_{X,c})$
and $\fN(\cI_{X_{\pi(c)},c})$ denote the diagrams of initial exponents with respect
to the coordinates $x$ and $v$ (respectively) for $Z$ and $Z_{\pi(c)}$ (respectively).
By semicontinuity of the diagram of initial exponents, we can assume that
$\fN(\cI_{X,c})$ and $\fN(\cI_{X_{\pi(c)},c})$ are constants, say
$\fN \subset \IN^n$ and $\fN^* \subset \IN^{n-p}$ (respectively), on the closed
points $c$ of $V$. It follows in a simple way
that $\fN
= \IN^p \times \fN^*$, and $\fN(\cI_{X_\eta,b})
= \fN^*$, where $\fN(\cI_{X_\eta,b})$ is the diagram with respect to the coordinates
$v = (v_1,\ldots,v_{n-p})$ for $N_\eta$. (Compare with \cite[Proof of Th.\,6.18]{BMda1}.)
In particular, $H_{X,a} = H_{X_\eta,b}^{(p)}$. ($p$ is the transcendence degree of $\um$
over $\uk$.)

A presentation $\ucI$ of the Hilbert-Samuel function $H_{X,\cdot}$ at $a$ with respect
to the coordinates $x=(u,v)$, as constructed in \cite[Ch.\,III]{BMinv}, is characterized
by certain formal properties \cite[(7.2)]{BMinv} related to the vertices of
$\fN(\cI_{X,c})$ above. Because of the product structure $\fN = \IN^p \times \fN^*$
of this diagram, it is easy to verify that, if these properties are satisfied
at every closed point of an open
subset of $W$, then they are satisfied by the induced marked ideal $\ucI_\eta$
at $b$. The details are left to the reader.
\hfill$\square$

\section{Absolute desingularization}
In this section, we apply the same approximation methods of \cite[\S8]{EGAIV}
that we used in the proof of Theorem 4.5 to show that any
desingularization algorithm for $\IQ$-varieties that is functorial with
respect to smooth morphisms extends
uniquely to a desingularization algorithm for a class $\fC$ of schemes over
$\IQ$ which includes all
varieties of characteristic zero as well as their localizations and
Henselizations along closed subvarieties. Moreover,
the algorithm for $\fC$ will be functorial with respect to all regular
morphisms between schemes in $\fC$.

We refer to \cite[\S\S18.6,\,18.8]{EGAIV} for definitions of
\emph{Henselization} and \emph{strict Henselization}.

\begin{remarks}\label{Crem1}
(1) One of our main motivations here is to extend the desingularization algorithm
for a variety $X$ to its Henselization $X_Z^h$ along a closed subvariety $Z$.
Since $X_Z^h\to X$ is a regular morphism, we could just
pull back the desingularization sequence from $X$, but it would not be clear
that the induced desingularization sequence depends only on the
scheme $X_Z^h$. The problem is that, while the ground field morphism
$X\to\Spec(\uk)$ is more or less unique (by \S4.1), the morphism
$X_Z^h\to\Spec(k)$ admits many deformations in general.

For example, even in the case $X=\IA^1_\IC$, if $x$ is the origin, then the
homomorphism $\cO_{X,x}\to
\ka(x)\xrightarrow{\cong}\IC$ admits many different extensions to the
Henselization.

Therefore, given a desingularization algorithm for $\uk$-varieties, the
blowing-up sequence for $X^h_Z$ obtained by pulling back that of $X$ might
depend on the morphism $X^h_Z\to\Spec(\uk)$.
We overcome this obstacle by descending to $\IQ$ --- we show that an absolute
desingularization algorithm for varieties defined over $\IQ$ induces a
desingularization algorithm for Henselian varieties (and certain
other schemes) that depends only on the schemes.
\smallskip

(2) Our Henselian result will be used in \cite{Tem2} to construct a canonical
desingularization of rig-regular formal varieties in characteristic zero
(independent of algebraization). The class includes, for
example, formal completion of a variety along its singular locus.

It seems to be an interesting open question whether the algorithm of \cite{BMinv}
extends to functorial desingularization of formal varieties in general. It is
true that, if $X$ and $Y$ are varieties (over perhaps different ground fields)
and $\wcO_{X,x} \cong \wcO_{Y,y}$, for some $x \in X,\, y \in Y$,
then the desingularizations of $X,\,Y$ induce the same sequences of formal
blowings-up of $\wX_x := \Spec \wcO_{X,x}$ and $\wY_y$.
We can show this using a formal presentation of the Hilbert-Samuel function
as given in \cite[\S78]{BMinv} together with the marked ideal techniques of
\cite{BMfunct} and commutativity of blowing up and formal completion
\cite[Lemma 2.1.8]{Tem}.
\end{remarks}

\begin{definition}\label{Cdef}
Consider a filtered projective family $\{X_i\}_{i\in I}$ of $\IQ$-varieties
with smooth affine transition morphisms $f_{ji}:X_j\to X_i$. The
projective limit $X=\projlim_{i\in I}X_i$ exists in the category of schemes, by
\cite[\S8]{EGAIV}. Assume that $X$ is Noetherian. Since the morphisms
$X_j\to X_i$, $j\ge i$, are regular, each projection $f_i:X\to X_i$ is regular.
(See, for example, \cite[\S1.4]{Sw}. The same argument shows, moreover, that
$X$ is Noetherian provided that $\dim X_i$ is bounded.) Let $\fC_\loc$ denote
the family of all such schemes $X$, and let
$\fC$ we denote the class of schemes each obtained by gluing together
finitely many elements of $\fC_\loc$.
\end{definition}

\begin{remarks}\label{Crem2}
(1) A simple argument in the proof of Theorem \ref{Cth2} below shows that we
could consider only families of affine varieties $X_i$ in Definition \ref{Cdef}
--- we would get a smaller category $\fC_\loc$, while the category
$\fC$ would not change.
\smallskip

(2) Any noetherian (or even quasi-compact quasi-separated) scheme in characteristic
zero is a projective limit of $\IQ$-varieties, by a noetherian approximation
theorem of Thomason \cite[C.9]{Th}, but the transition
morphisms are not smooth (or even flat) in general.

For example (in positive characteristic), if $K$ is a perfect field of
positive transcendence degree over $\IF_p$, then $\Spec(K)$
is not a projective limit of a filtered family of $\IF_p$-varieties with
smooth transition morphisms. This follows from the fact that $K$ is not
separable over any finitely generated subfield of positive
absolute transcendence degree.

A deep theorem of Popescu \cite{Pop} states that any regular morphism
$X\to Y$ is a projective limit of smooth morphisms $X_i\to Y$. But the
transition morphisms $X_j\to X_i$ and the projections $X\to X_i$ cannot be
made regular in general.
\end{remarks}

\begin{theorem}\label{Cth}
Any desingularization algorithm for $\IQ$-varieties that is functorial with respect to smooth morphisms extends
uniquely to a desingularization algorithm on $\fC$ that is functorial with respect to all regular morphisms. 
Moreover, if the original algorithm satisfies the stronger conditions of 
Theorem 6.1, then the extended
algorithm satisfies the same conditions.
\end{theorem}

\begin{proof}
Fix a desingularization algorithm $\cF$ for $\IQ$-varieties. First we extend $\cF$ to $\fC_\loc$. Let $X$ be an
element of $\fC_\loc$ and let $X=\projlim_{i\in I} X_i$ denote a representation of $X$ as a projective limit of
$\IQ$-varieties with smooth affine transition morphisms. Then the desingularizations $\cF(X_i)$ are compatible,
so that each of them induces the same desingularization sequence for $X$, which we denote $\cF(X)$. 
Moreover, if the $\cF(X_i)$  satisfy the conditions of Theorem 6.1, then 
$\cF(X)$ also satisfies them.

We have to
prove that $\cF(X)$ is independent of the choice of the projective limit
representation and that this extension of
$\cF$ to $\fC_\loc$ is compatible with all regular morphisms. For both tasks,
it is enough to prove that, given another family $\{Y_j\}_{j\in J}$ of
$\IQ$-varieties with smooth affine transition
morphisms and limit $Y$, and given a regular morphism $h:Y\to X$,
there exist $i\in I$, $j\in J$ and a regular
morphism $h_{ji}:Y_j\to X_i$ compatible with $h$ in the sense that the
following diagram commutes.
$$
\begin{CD}
Y  @>>> Y_j \\
@V h VV    @VV h_{ji} V \\
X  @>>> X_i\\
\end{CD}
$$
Indeed, if we have morphisms such that the diagram commutes, then, on the one
hand, $\cF(Y)$ is induced by $\cF(Y_j)$ and hence
by $\cF(X_i)$ (since $\cF$ is compatible with the smooth morphisms $h_{ji}$),
and, on the other hand, $\cF(X)$ is induced by $\cF(X_i)$. Therefore,
$\cF(Y)$ is induced by $\cF(X)$, thus proving compatibility with regular
morphisms. The fact that $\cF(X)$ is well defined
is then obtained by applying the preceding argument to an isomorphism
$X\xrightarrow{\cong} X$ and two representations of $X$ as a projective
limit.

To find a regular $h_{ji}$ as above:  Fix $i$. By \cite[Cor.\,8.13.2]{EGAIV}, the morphism $Y\to X_i$ factors
through $Y_j$, for some $j$. So we get a morphism $h_{ji}$ and it remains to show 
only that it can be chosen regular. We
claim that the image of $Y$ in $Y_j$ lies in the maximal open subscheme $U$ of
$Y_j$ such that $h_{ji}|_U$ is smooth.
Indeed, let $y\in Y$ and
let $y_j\in Y_j$, $x_i\in X_i$ denote its images. Then the local homomorphisms 
$\psi: \cO_{Y_j,y_j} \to \cO_{Y,y}$
and $\cO_{X_i,x_i} \to \cO_{Y,y}$ are regular (since the morphisms $Y\to Y_j$ 
and $Y\to X\to X_i$ are regular);
hence the homomorphism $\cO_{X_i,x_i} \to \cO_{Y_j,y_j}$ is regular,
by \cite[Lemma 33.B]{Mat} (where all we
need to know about $\psi$ is that it is faithfully flat). Thus $y_j\in U$ and 
therefore $Y$ lies
in $U$. Applying \cite[Cor.\,8.13.2]{EGAIV} again, we see that the morphism $Y\to U$ factors through $Y_k$
for some $k\ge j$. Then the morphism $h_{ki}:Y_k\to U\to X_i$ is smooth since $Y_k\to Y_j$ and $U\to X_i$ are
smooth.

Finally, we can extend $\cF$ from $\fC_\loc$ to $\fC$ using the gluing
argument of \cite[Prop.\,3.37]{Ko}:
Given $X$ in $\fC$, take a covering of $X$ by open subschemes
$X_i \in \fC_\loc$, $i=1,\ldots k$. Then the disjoint unions $\coprod_i X_i$
and $\coprod_{i\leq j}X_i\cap X_j$ belong to $\fC_\loc$, and we have a
commutative diagram
$$
\begin{CD}
\coprod_{i\leq j}X_i\cap X_j @>>> \coprod_i X_i\\
@VVV                              @VVV\\
\coprod_i X_i     @>>>            X\\
\end{CD}
$$
where the left and top arrows are induced by $X_i\cap X_j \hookrightarrow X_i$
and $X_i\cap X_j \hookrightarrow X_j$, respectively. These arrows are \'etale,
hence smooth, so that $\cF$ on $\fC_\loc$ is compatible with them. It follows
that the blowing-up sequences $\cF(X_i)$ glue together to give a
desingularization $\cF(X)$. Clearly,
$\cF$ is compatible with regular morphisms among members of $\fC$.
\end{proof}

\begin{theorem}\label{Cth2}
(1) The class $\fC_\loc$ contains all
separated varieties of characteristic zero,
as well as their localizations,
Henselizations and strict Henselizations along closed subvarieties.
As a result, the class $\fC$ contains analogous classes of schemes
(that are not necessarily separated).

(2) If $\{X_i\}_{i\in I}$ is a filtered projective family of separated schemes in
$\fC_\loc$ with regular
affine transition morphisms then $X=\projlim_{i\in I}X_i$ belongs to $\fC_\loc$.
\end{theorem}

\begin{proof}
We start with (2). For each $X_i$ fix a representation
$X_i\xrightarrow{\cong}\projlim_{j\in J_i}X_{ij}$ with smooth affine
transition morphisms between $\IQ$-varieties $X_{ij}$. Using \cite[C.7]{Th} we can assume that all $X_{ij}$
are separated. We recall that the \emph{schematic image} $X'_{ij}$ of $X_i$ in
$X_{ij}$ means the smallest closed subscheme of
$X_{ij}$ through which the morphism $X_i\to X_{ij}$ factors.
The morphism $X_i\to X'_{ij}$ is regular and each
transition morphism $X_{ij}\to X_{ik}$ restricts to a regular morphism
$X'_{ij}\to X'_{ik}$; hence, by replacing
$X_{ij}$ with $X'_{ij}$, for all $j$, we can assume that the projections
$X_i\to X_{ij}$ are schematically dominant.

We will now add certain transition morphisms $X_{ij}\to X_{kl}$ with $i\geq k$,
which are regular, affine, and moreover
make the entire family $\{X_{ij}\}_{i\in I,j\in J_i}$ into a filtered family with projective
limit $X$. This
will prove (2). For each $i'\geq i$ and $j\in J_i$ let $f_{i'ij}:X_{i'}\to X_{ij}$
denote the morphism obtained by composing
$X_{i'}\to X_i$ and $X_i\to X_{ij}$. Note that if $f_{i'ij}$ factors through a morphism
$f_{i'j'ij}:X_{i'j'}\to X_{ij}$, then $f_{i'j'ij}$ is unique because $X_{i'}\to X_{i'j'}$ is schematically
dominant and $X_{ij}$ is separated. If such $f_{i'j'ij}$ exists and is affine and regular, then we declare that
$(i'j')\ge(ij)$. Affineness and regularity are preserved by composition,
so this defines an order on the set
$J:=\coprod_{i\in I}J_i$. Moreover, this makes $J$ into a filtered ordered set because the argument from the
proof of Theorem \ref{Cth} shows that for each $i'\ge i$ and $j\in J_i$ the morphism $f_{i'j'ij}$ exists and is
regular and affine for $j'\geq j'_0(i,i',j)$. Since
$X\xrightarrow{\cong}\projlim_{i\in I}\projlim_{j\in
J_i}X_{ij}\xrightarrow{\cong}\projlim_{(ij)\in J}X_{ij}$,
the family $\{X_{ij}\}_{(ij)\in J}$ is as required, and $X$ is in
$\fC_\loc$.

The assertion (1) follows from (2): Indeed, suppose that $Y$ is a variety
over a field $\ul$ that
is finitely generated over $\IQ$. Then $Y$ is a pro-open
subscheme of a $\IQ$-variety; hence $Y$ is the projective limit of all its open neighborhoods and the
transition morphisms are open immersions. So $Y$ is in $\fC_\loc$.
An arbitrary variety $X$ is of the
form $Y\otimes_{\ul} \uk := Y \times_{\Spec \ul} \Spec \uk$ with
$Y$ and $\ul$ as above (see Theorem 3.1),
so $X$ is the projective limit of varieties $X_i=Y\otimes_{\ul} \uk_i$
where $\uk_i$ is a finitely generated $\ul$-subfield of $\uk$. The transition morphisms are regular by the
characteristic zero assumption, and if $X$ is separated then each $X_i$ is a separated $\uk_i$-variety. Then all
$X_i \in \fC_\loc$ and hence $X \in \fC_\loc$, by (2). Finally, the strict Henselization
(respectively, Henselization, or localization) of $X$ along a closed subvariety
$Z$ is a projective limit of a family of $X$-\'etale
schemes $X_j$. Since the $X_j$ are separated $\uk$-varieties, we get (1) using (2)
again.
\end{proof}

\begin{remarks}\label{Crem3}
(1) It is interesting to ask whether the category $\fC$ can be naturally extended further. (See, for example, Remark \ref{Crem1}(2).)
\smallskip

(2) In principle, if $X$ admits a regular morphism $f$ to a variety $Y$,
we could induce a desingularization of $X$ from
a desingularization of $Y$ (even though $X$ might not be quasi-excellent!). We do not know if this desingularizaton would be independent of $f$; even a tool as strong as Popescu's theorem would seem to be of no help here.
(See also remark \ref{Crem2}(2).)
\end{remarks}

\end{document}